\documentclass{amsart}
\usepackage{latexsym, amsfonts, euscript}
\usepackage{epsfig}
\usepackage{graphicx,color}
\usepackage{pst-plot}
\newtheorem{theorem}{Theorem}[section]
\newtheorem{proposition}[theorem]{Proposition}
\newtheorem{definition}[theorem]{Definition}
\newtheorem{lemma}[theorem]{Lemma}
\newtheorem{corollary}[theorem]{Corollary}

\newtheorem{remark}[theorem]{Remark}

\newtheorem{example}[theorem]{Example}
\newtheorem{examples}[theorem]{Examples}

\def\N{{\mathbb N}}

\def\R{{\mathbb R}}

\begin{document}
\title{A Geometric Approach to Saddle Points of Surfaces}

\author{Sudhir R. Ghorpade}
\address{Department of Mathematics, % \newline  \indent
Indian Institute of Technology Bombay,\newline \indent
Powai, Mumbai 400076, India.}
\email{srg@math.iitb.ac.in}
\urladdr{http://www.math.iitb.ac.in/$\sim$srg/}

\author{Balmohan V. Limaye}
\address{Department of Mathematics, % \newline  \indent
Indian Institute of Technology Bombay,\newline \indent
Powai, Mumbai 400076, India.}
\email{bvl@math.iitb.ac.in}
\urladdr{http://www.math.iitb.ac.in/$\sim$bvl/}

%\begin{abstract}
%We outline an alternative approach to the geometric notion of a saddle point for real-valued functions of %several 
%two variables. It 
%%is based on a geometric viewpoint and it 
%is argued that this is more natural compared to the usual treatment of this topic %based on Calculus alone. 
%in standard texts on Calculus. 
%%We also discuss the related problem of characterizing positive definite 
%%as well as nonnegative definite quadratic forms. For quadratic forms in $3$ variables, we give an elementary and self-contained proof of Sylvester's %Minorant 
%%Criterion for positive definiteness as well as for 
%%nonnegative definiteness. 
%\end{abstract}

\maketitle
 
\section{Introduction}
\label{sec:intro}

What is a saddle point of a surface in $3$-space? 
A reasonable answer is: %to this question can be as follows. %either of the following. A 
a saddle point is like %the tip of 
the center point of a horse saddle or the low point of a ridge joining two peaks. 
%If we think of a surface in $3$-space or (the graph of) a real valued function of two real variables, then 
In other words, a saddle point is that  peculiar point on the surface which is at once a peak along a path on the surface and a dip along another path on the surface. Another answer that is %more 
mundane but %one that is 
more likely to fetch points in a Calculus test %and it goes 
is as follows. A \emph{saddle point} of a real-valued function of two real variables 
is %usually defined  (cf. \cite[\S 9.9]{Ap} or \cite[\S 3.3]{MT}) as 
a critical point (that is, a point where the gradient vanishes) which is not a local extremum. 
%The two answers seem to complement each other. Indeed, the first tells us what a saddle point should be and the second %tells us what it is. In other words, the 
The first answer gives an intuitive description of a saddle point, while the second is the mathematical definition commonly given in most texts on Calculus. %(See, for example, \cite[\S 9.9]{Ap} or \cite[\S 3.3]{MT}.) 
(See, e.g., \cite[\S 9.9]{Ap} or \cite[\S 3.3]{MT}.) 
A typical example %to illustrate both the answers is that of 
is the hyperbolic  paraboloid given by $z=xy$ or by (the graph of) the function $f:\R^2\to \R$ defined by $f(x,y):=xy$. Here the origin is a saddle point. Indeed if we look at the paths along the diagonal lines $y=-x$ and $y=x$ in the plane, 
% and the corresponding paths on the surface  $z=xy$ 
then we readily notice that the origin is at once a peak and a dip. Also, the origin is the only critical point of $f$ and clearly $f$ does not have a local extremum at the origin. 

The aim of this article is first, to point out that there is a significant disparity between the two answers, and second, to suggest an 
alternative approach to saddle points which may take care of this. 
The first point is easy to illustrate. There are surfaces or rather, functions of two variables where the conditions in the second answer are met but the geometric picture is nowhere close to the description in the first answer. For 
example, if $f:\R^2\to \R$ is defined by $f(x,y):=x^3$ or by $f(x,y):=x^2+y^3$, then the origin is a saddle point according to the usual mathematical definition, but the corresponding surface (Figure \ref{fig:x2plusy3}) hardly looks like a 
saddle that you might want to put on a horse for any rider! 
Another unsatisfactory aspect %of the usual definition of a saddle point 
is the {\it a priori} assumption that the saddle point is a critical point, 
that is, a point at which the gradient exists and is zero. % at that point. 
This is quite unlike the usual definitions of analogous concepts in one variable calculus, such as local extrema or points of inflection, where one makes a clear distinction between a geometric concept and its %the 
analytic characterization 
%when conditions such as differentiability hold. 
(See, for instance, \cite{GL} and its review \cite{W}.)
%%percented
%Apart from such anomalies, there is also an aspect of the usual definition of a saddle point that is somewhat unsatisfactory. To see this, let us first note that most texts and teachers of Calculus would take pains to point out the distinction between a geometric concept and the analytic characterization or criterion employed in studying them. (See, for instance, \cite{GL} and its review \cite{W}.) For example, a real-valued function $f$ of one variable defined in the neighborhood of a point $c$ is said to have a {\em local extremum} at $c$ if $f(x)$ is bounded above or bounded below by $f(c)$ for all $x$ in an open interval around $c$. 
%%Likewise, we would say that $f$ has a {\em point of inflection} at $c$ if $f$ changes concavity near $c$. 
%In these cases we don't need in the definition that $f$ is differentiable %or twice differentiable 
%at $c$,  but when it is so we obtain as a consequence that $f'(c)$ %or $f''(c)$ 
%must, in fact, be zero. However, when it comes to saddle points, we begin by assuming the point to be a critical point, that is,  the gradient should exist and vanish at that point. In this sense the usual definition is not geometric and is a departure from the good practices adopted in the theory of functions of one or more real variables. 
The definition we propose here seems to fare better on these counts in the case of functions of two variables. The basic idea is quite simple and, we expect, scarcely novel. However, we have not seen in the literature an exposition along the lines given here. For this reason, and with the hope that the  treatment 
suggested here could become standard, we provide a fairly detailed discussion of the definition, basic results and a number of examples in the next three sections.  Alternative approaches and extensions are briefly indicated in a remark at the end of the paper and we thank the referee for some of the suggestions therein. 

\section{Definition of a Saddle Point}

%Let us begin with some preliminary notions. 
Let $D$ be a subset of $\R^2$. %$D \subseteq \R^2$. 
A \emph{path} in $D$ is a continuous map from  $[a,b]$ to $D$. Here, and hereafter, while writing open or closed intervals such as $(a,b)$ or $[a,b]$, it is tacitly assumed that $a,b\in \R$ with $a<b$. Given any $\mathbf{p}\in D$, 
a path $\gamma : [a,b]\to D$ is said to \emph{pass through} $\mathbf{p}$ if $\gamma(t_0) = \mathbf{p}$ for some $t_0\in (a,b)$. 
A path  $\gamma : [a,b]\to D$  is said to be \emph{regular} if $\gamma$ is differentiable on $(a,b)$ and $\gamma'(t) \ne \mathbf{0}$ for all $t\in (a,b)$. 
Two regular paths $\gamma_1 : [a_1,b_1]\to D$ and $\gamma_2 : [a_2,b_2]\to D$ are said to \emph{intersect transversally} at some $\mathbf{p}\in D$ if there are $t_i\in (a_i, b_i)$ such that $\gamma_i(t_i) = \mathbf{p}$ for $i=1,2$ and moreover, 
$\gamma_1'(t_1)$ and $\gamma_2'(t_2)$ are not multiples of each other. 
In other words, the two paths pass through $\mathbf{p}$ and their tangent vectors at $\mathbf{p}$ are not parallel. 

\begin{examples}
\label{Exa:2.1}
{\rm
(i) $\gamma : [-1,1]\to \R^2$ defined by $\gamma(t):= (t, t^2)$ is a regular path, while 
$\widetilde{\gamma }: [-1,1]\to \R^2$ defined by $\widetilde{\gamma}(t):= (t^2, t^3)$ is not a regular path. 

(ii) If $\gamma_1, \gamma_2 : [-1,1]\to \R^2$ are defined by 
$\gamma_1(t):= (t, -t)$ and $\gamma_2(t):= (t, t)$, then 
$\gamma_1$ and $\gamma_2$ are regular paths in $\R^2$ which intersect transversally at the origin. Further, the path $\gamma_3:[-1,1]\to \R^2$ defined by $\gamma_3(t):= (2t+t^2, 2t - t^2)$, %then $\gamma_3$ 
is also regular and passes through the origin. The paths $\gamma_1$ and $\gamma_3$ intersect transversally at the origin, whereas the paths $\gamma_2$ and $\gamma_3$ do not. % intersect transversally at the origin.
% On the other hand,  
}
\end{examples}

\begin{figure}
	\centering
	 {\scalebox{0.52}{\includegraphics{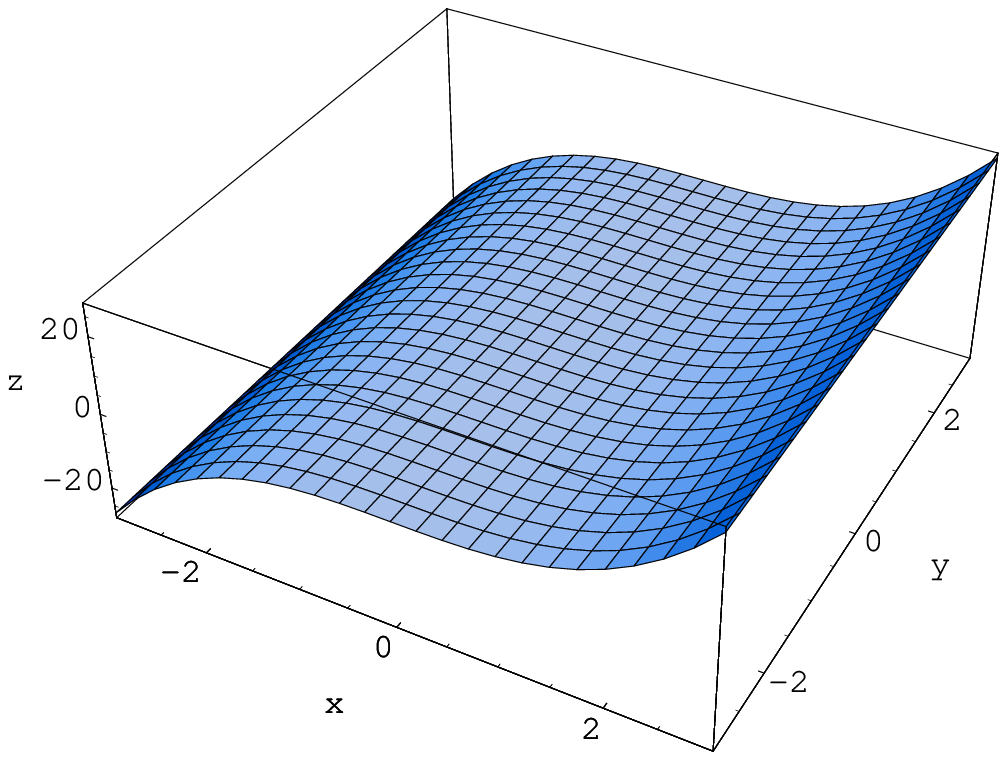}}}
		\qquad 
		{\scalebox{0.52}{\includegraphics{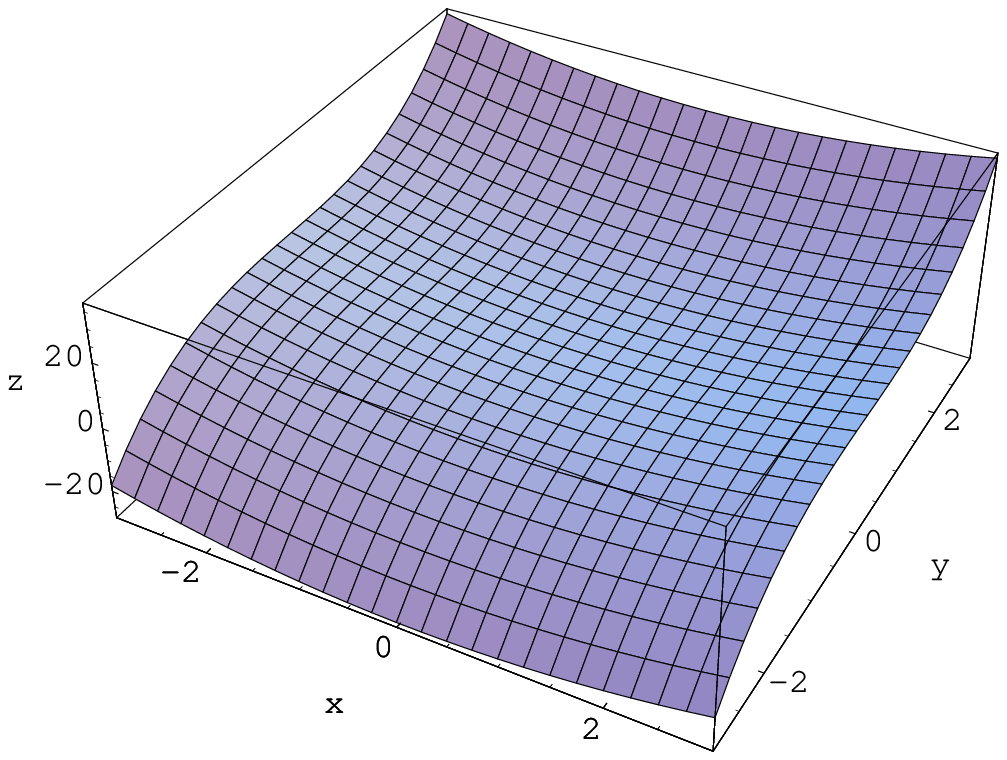}}}
	\caption{Graphs of $z=x^3$ and $z= x^2 + y^3$ near the origin}
	\label{fig:x2plusy3}
%	Illustration of a fake saddle point: $f(x,y) = x^2 + y^3$ near the origin}
\end{figure}

Let $D \subseteq \R^2$, $\mathbf{p}\in D$ and %$\mathbf{p}$ be an interior point of $D$.
%A path Let 
$\gamma : [a,b]\to D$ be a regular path in $D$ passing through $\mathbf{p}$ so that %. Thus %$\mathbf{p}\in D$ and 
$\gamma(t_0) = \mathbf{p}$ for some $t_0\in (a,b)$. 
%Notice that since $\gamma$ is regular, %it is one-one on $(a,b)$ and 
%the point $t_0$ is uniquely determined by $\mathbf{p}$. 
Now, any 
$f:D\to \R$ can be restricted to (the image of) $\gamma$ so as to obtain 
a real-valued function of one variable $\phi:[a,b]\to \R$ defined by
$\phi (t):= f(\gamma (t))$. We shall say that $f$ has a local maximum (resp: local minimum) at $\mathbf{p}$ along $\gamma$ if $\phi$ has a local maximum (resp: local minimum) at $t_0$. 

\begin{definition}
\label{Def:2.2}
{\rm
Let $D \subseteq \R^2$ and $\mathbf{p}$ %\in D$. 
 be an interior point of $D$. 
A real-valued function 
$f:D\to \R$ has %is said to have 
a \emph{saddle point} at $\mathbf{p}$ if there are regular paths $\gamma_1$ and $\gamma_2$ in $D$ intersecting transversally at $\mathbf{p}$ such that $f$ has a local maximum at $\mathbf{p}$ along $\gamma_1$, 
while $f$ has a local minimum at $\mathbf{p}$ along $\gamma_2$. 
}
\end{definition}

The above definition is a faithful abstraction of the idea that a saddle point is the point at which the graph of the function is at once 
a peak along a path and a dip along another path. The condition that the two paths intersect transversally might seem %rather 
technical. But its significance will be clear from Example \ref{Exa:2.3}(iii) below. 
%a moment's reflection should convince the reader of its importance. 
%%%%%%%%%%%%%%%%%%%%%%%%%%%%%%%%%%%%%%%
%A moment's reflection should also 
%convince the reader the importance of the condition that 
%%it is important to require that 
%the two paths along which the point becomes a peak or a valley %ought to 
%intersect transversally. 
%%%%%%%%%%%%%%%%%%%%%%%%%%%%%%%%%%%%%%%%%
%For instance, this condition ensures %also 
%that there are no saddle points when $n=1$ because on the real line, no two regular paths intersect transversally. Example \ref{Exa:2.3}(iii) below illustrates the significance of this condition when $n=2$.

It may be remarked that in our definition of a saddle point, we have permitted ourselves as much laxity as is usual while defining local extrema. % of functions of two variables. 
To wit, %a point around which a 
if a function is locally constant  at $\mathbf{p}$, then it has a local 
maximum as well as a local minimum at $\mathbf{p}$. In the same vein, a 
locally constant function at $\mathbf{p}$ has a saddle point at  $\mathbf{p}$. More generally, if a function is locally constant along 
two regular paths intersecting transversally at $\mathbf{p}$, then it has a 
saddle point at $\mathbf{p}$. If we don't want to be so indulgent, then we  can use the stronger notion of a \emph{strict} saddle point. A  \emph{strict saddle point} is defined simply  by replacing in Definition  \ref{Def:2.2}, 
%local extrema by strict local extrema. 
local maximum by strict local maximum and local minimum by strict local minimum. Indeed, it is the notion of a strict saddle point that comes closest to one's geometric intuition about saddle points. In almost all the examples as well as the criteria for saddle points discussed here, it is %will be 
seen that the function has, in fact, a strict saddle point. 

\begin{figure}
	\centering
	 {\scalebox{0.52}{\includegraphics{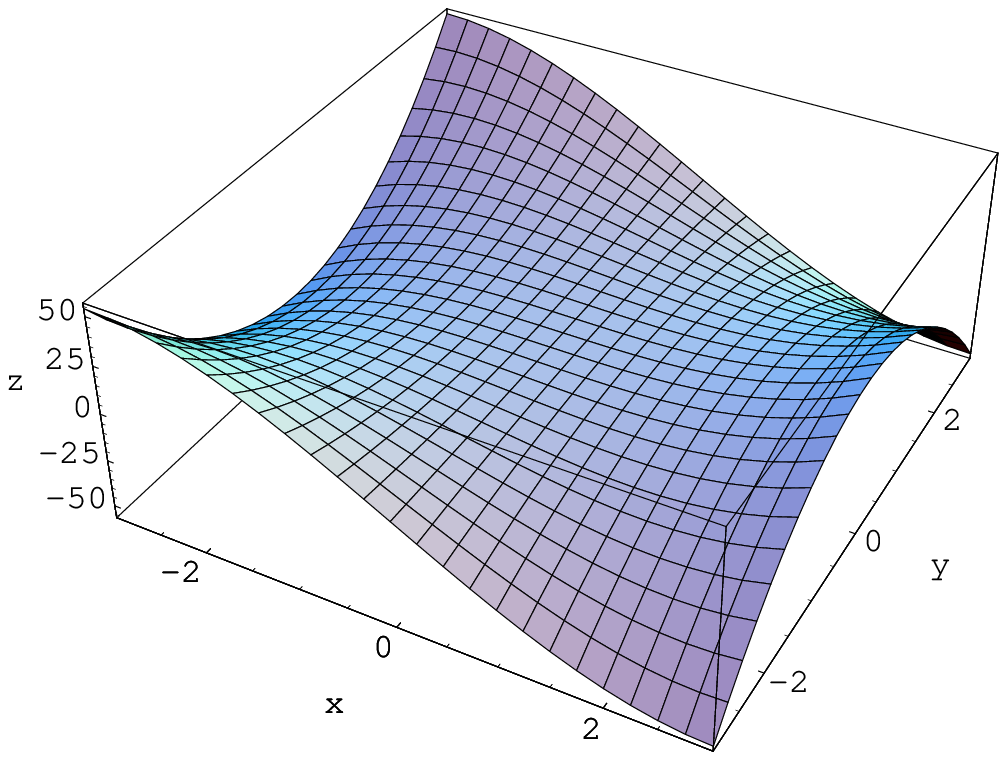}}}
		\qquad 
		{\scalebox{0.52}{\includegraphics{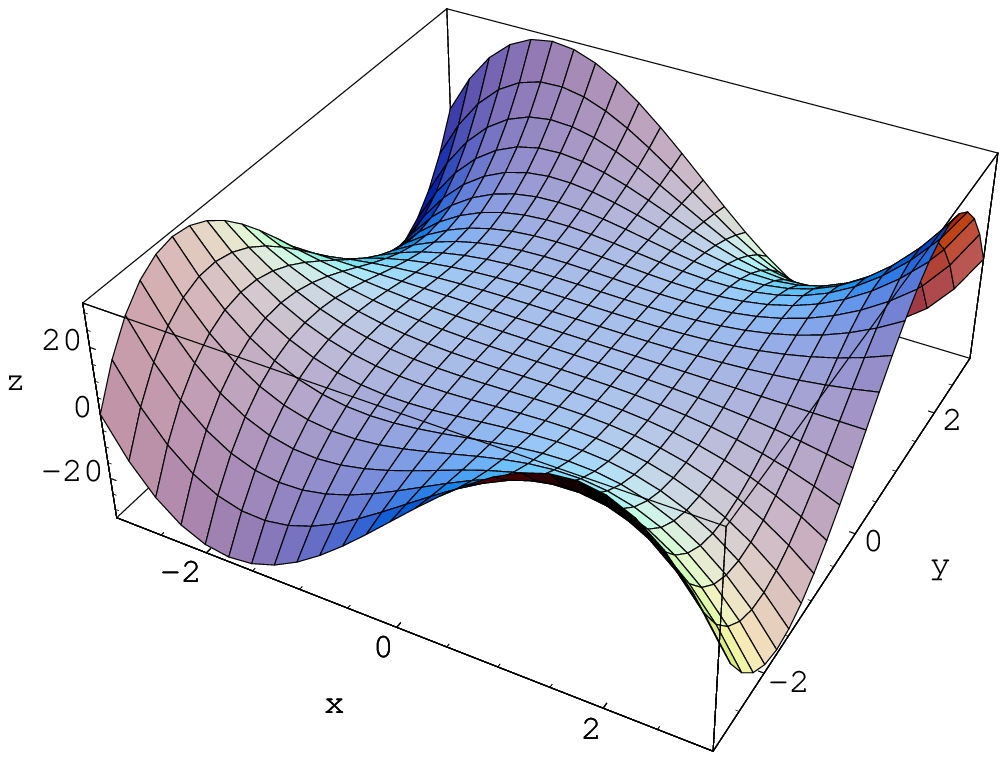}}}
	\caption{Monkey saddle $z=x^3-3xy^2$ and dog saddle $z= x^3y - xy^3$}
	\label{fig:twosaddles}
%	Illustration of a fake saddle point: $f(x,y) = x^2 + y^3$ near the origin}
\end{figure}

\begin{examples}
\label{Exa:2.3}
{\rm
(i)  [Hyperbolic paraboloid] %Consider 
The function $f:\R^2\to \R$ defined by 
$f(x,y):=xy$ has a saddle point at $(0,0)$. %the origin. 
To see this, it suffices to consider the paths $\gamma_1$ and $\gamma_2$ in Example \ref{Exa:2.1}(ii). Similarly, one can see % argument shows 
that if $a,b,c,d\in \R$ with $ad-bc\ne 0$, then $f:\R^2 \to \R$ defined by $f(x,y):=(ax+by)(cx+dy)$ has a saddle 
point --- in fact, a strict saddle point, at $(0,0)$. 

(ii) [Monkey saddle]
%%%%%%%%%%%%%%% Older Example: Without the 3 %%%%%%%%%%%%%%%%%%%%%%%%%%%% 
%The function $f:\R^2\to \R$ defined by 
%$f(x,y):=x^3 - xy^2$ has a saddle point at the origin. This can also be seen rather trivially, by considering the paths $\gamma_1$ and $\gamma_2$ in Example \ref{Exa:2.1}(ii). In fact, $f$ vanishes along these two paths and hence has a local maximum as well as a local minimum at the origin. It is natural to ask if one can come up with nontrivial paths and whether $f$ has a strict saddle point at the origin. To this end, it helps to look at the level curves of $f$. These suggest that the parabolic paths given by $t\mapsto (t^2,\,  t)$ and 
%$t\mapsto (t - t^2, \, t+t^2)$ should do the needful. %We leave it to the reader to 
%Indeed, it can be readily checked that for these paths the conditions in 
%Definition \ref{Def:2.2} are satisfied. Moreover, $f$ has a strict 
%saddle point at the origin. 
%%%%%%%%%%%%%%%%%%%%%%%%%%%%%%%%%%%%%%%%%%%%%%%%%%%%%%%%%%%%%%%%%%%%
The function $f:\R^2\to \R$ defined by 
$f(x,y):=x^3 - 3xy^2$ has a strict saddle point at the origin.  To prove this, it helps to look at the level curves of $f$. We then find that it suffices to consider the parabolic paths given by $t\mapsto (- t \sqrt{3} + t^2,\;  t + t^2 \sqrt{3})$ and $t\mapsto ( t \sqrt{3} - t^2, \; t+ t^2 \sqrt{3})$ for $t\in [-\sqrt{3}, \sqrt{3}]$.
% should do the needful. 
The surface $z=f(x,y)$ or the graph of $f$ near the origin is shown in Figure \ref{fig:twosaddles} on the left. It may be interesting to try and visualize these paths on this surface. 
% %We leave it to the reader to 
%Indeed, it can be readily checked that for these paths the conditions in 
%Definition \ref{Def:2.2} are satisfied. Moreover, $f$ has a strict 
%saddle point at the origin. These results may be compared with the graph of $f$ near the origin, which is shown in Figure \ref{fig:twosaddles} on the left. 
%%this works. %everything works out. 
%%$\gamma_1$ and $\gamma_2$ in Example \ref{Exa:2.1}(i). 

(iii) [Fake saddle] Consider $f:\R^2\to \R$ defined by $f(x,y):=x^3$. 
%The paths 
In this case $\gamma_1,\gamma_2: [-1,1]\to \R$ defined by $\gamma_1(t):= (-t^2, t)$ and $\gamma_2(t)=(t^2,t)$ are regular paths passing through $\mathbf{0}$. Also, %the function 
$f$ has a strict local maximum at $\mathbf{0}$ along $\gamma_1$ and a strict local minimum at $\mathbf{0}$ along $\gamma_2$. However, $\gamma_1$ and $\gamma_2$ do not intersect transversally at $\mathbf{0}$. In fact, as the surface on the left in Figure \ref{fig:x2plusy3} indicates, $f$ does not have a saddle point at $\mathbf{0}$. A formal proof of this 
%in a slightly more general case 
is given later in Example \ref{Exa:4.1}(iii). 
}
\end{examples}

%As alluded to in the Introduction, our definition of a saddle point does not presuppose that 
%the point is a critical point. However, in the case of differentiable functions of two variables, 
We now show that a saddle point is necessarily a critical point. 
%%%%%%% CHANGE: The line below percented.
%For functions of $n$ variables, a weaker result holds. 
In what follows, by $\nabla f (\mathbf{p})$ we denote 
the gradient of a function $f$ %of $n$ variables 
at an interior point $\mathbf{p}$ of its domain.  

%%%%%%% CHANGE: Revised Prop and an added Remark thereafter. 
\begin{proposition}
\label{Prop:2.4}
Let $D \subseteq \R^2$ and $\mathbf{p}$  be an interior point of $D$. 
If $f:D\to \R$ is differentiable at $\mathbf{p}$ and has a saddle point  at $\mathbf{p}$, then  
$\nabla f (\mathbf{p}) = \mathbf{0}$. 
\end{proposition}

\begin{proof}
For $i=1,2$, let $\gamma_i : [a_i,b_i]\to D$ satisfy the conditions in 
Definition \ref{Def:2.2} with $n=2$, and let us write $\mathbf{p} = \gamma_i(t_i)$ with $t_i\in (a_i, b_i)$ and $\phi_i:=f\circ \gamma_i$. Since $f$ is differentiable at $\mathbf{p}$ and $\gamma_i$ is regular, by the chain rule, $\phi_i'(t_i)$ exists and equals $\nabla f (\mathbf{p}) \cdot \gamma_i'(t_i)$ for $i=1,2$. On the other hand, since $\phi_i$ have local extrema at $t_i$, we have $\phi_i'(t_i) =0$ for $i=1,2$. 
%Further, since $\gamma_1$ and $\gamma_2$ are regular paths intersecting %transversally at $\mathbf{p}$,
Now, since %the tangent vectors 
$\gamma_1'(t_1)$ and $\gamma_2'(t_2)$ are linearly independent vectors in $\R^2$, we can conclude that %  It follows that 
$\nabla f (\mathbf{p}) = \mathbf{0}$. 
%span a plane $\mathcal{P}$ (say) in $\R^n$. 
%It follows that $\nabla f (\mathbf{p})$ is orthogonal to $\mathcal{P}$.  
\end{proof}

\section{Discriminant Test}

The Discriminant Test or the Second Derivative Test is a high point of any exposition of local extrema and saddle points of functions of two real variables. It facilitates easy checking of saddle points in many, but not all, cases. The classical definition of a saddle point given in the Introduction is, in fact, tailor-made so that the Discriminant Test can be proved easily. Some texts (e.g., \cite[p.347]{CJ}) even take an easier option to \emph{define} a saddle point as a critical point where the `discriminant' is negative. This may appear a bit like putting the cart before the horse. 
%Nonetheless, 
But the importance of the Discriminant Test can hardly be overemphasized and 
it seems imperative that it remains %continues to be 
available with our geometric notion of a saddle point. 

%Thus we prove in this section the Discriminant Test for functions of two variables.
%, and an analogous result for functions of $n$ variables. 
Let us recall that a {\em binary quadratic form} (over $\R$) is a polynomial of the form
$$
Q(\mathbf{h})=Q(h_1, h_2):= 
%\left[ \begin{array}{cc} h_1 & h_2 \end{array} \right]\left[ \begin{array}{cc} a & b \\b & c \end{array} \right]\left[ \begin{array}{c} h_1 \\ h_2 \end{array} \right] = 
ah_1^2 + 2bh_1h_2 + ch_2^2,
$$
where $\mathbf{h}=(h_1,h_2)$ is a pair of variables and $a,b,c$ are (real) constants. 
%The corresponding $2\times 2$ real symmetric matrix 
%$$
%A = \left[ \begin{array}{cc} a & b \\b & c \end{array} \right]
%$$
%is called the {\em matrix associated to} $Q$. Note that $Q(\mathbf{h})= \mathbf{h}A\mathbf{h}^T$. 
We say that $Q$ is \emph{positive definite} (resp: \emph{negative definite}) if 
$Q(\mathbf{u})>0$ (resp: $Q(\mathbf{u}) < 0$) for all $\mathbf{u}\in \R^2$, $\mathbf{u}\ne \mathbf{0}$. In case $Q$ takes positive as well as 
negative values, that is, if there are $\mathbf{u}, \mathbf{v}\in \R^2$ such that %$Q(\mathbf{u})>0$ and 
$Q(\mathbf{u})Q(\mathbf{v})<0$, then $Q$ is said to be \emph{indefinite}. In this situation, the vectors 
$\mathbf{u}$ and $\mathbf{v}$ are necessarily nonzero and they can not 
be multiples of each other since  %because
%$Q$ being homogeneous of degree $2$, one has
$Q(t \mathbf{h})= t^2 Q(\mathbf{h})$ for any $t\in \R$ and 
$\mathbf{h}\in \R^2$.

%%%%%%%%%%%%%%% A Lemma I had written earlier %%%%%%%%%%%%%%%%%%%%%%%%%
%\begin{lemma}
%\label{Lem:3.1}
%If $Q$ is an indefinite quadratic form in $n$ variables, then there are regular paths $\gamma_1, \gamma_2: [-1,1]\to \R^n$ with 
%$\gamma_1(\mathbf{0})=0= \gamma_2(\mathbf{u})$, such that 
%$Q\left(\gamma_1(t)\right) < 0$ and $Q\left(\gamma_2(t)\right) > 0$ for 
%all $t\in [-1,1]$ with $t\ne 0$. Moreover, $\gamma_1'(t)$ and $\gamma_2'(t)$ are not multiples of each other for any $t\in [-1,1]$. 
%\end{lemma}
%
%\begin{proof}
%Since $Q$ is indefinite, there are nonzero %vectors 
%$\mathbf{u}=(x_1,\dots , x_n)$ and $\mathbf{v}=(y_1, \dots , y_n)$ 
%in $\R^n$ such that $Q(\mathbf{u})>0$ and $Q(\mathbf{v})<0$. Define 
%$\gamma_1, \gamma_2: [-1,1]\to \R^n$ by 
%$\gamma_1(t):=(x_1t,\dots , x_nt)$ and $\gamma_2(t):=(y_1t,\dots , y_nt)$. Clearly, $\gamma_1$ and $\gamma_2$ are regular. %moreover, 
%Since $Q$ is homogeneous of degree $2$, for any $t\in [-1,1]$ with $t\ne 0$, we have 
%$$
%Q\left(\gamma_1(t)\right) = t^2 Q(\mathbf{u})< 0 \quad \mbox{and} \quad Q\left(\gamma_2(t)\right)= t^2 Q(\mathbf{v}) > 0. 
%%\quad \mbox{for all $t\in [-1,1]$ with $t\ne 0$.}
%$$
%%and moreover,
%In particular, 
%$\gamma_1'(t) = \mathbf{u}$ and %is not a multiple of
%$\gamma_2'(t) =\mathbf{v}$ are not multiples of each other 
%for any $t\in [-1,1]$.  
%\end{proof}
%%%%%%%%%%%%%%%%%%%%%%%%%%%%%%%%%%%%%%%%%%%%%%%%%%%%%%%%%%%%%%%%%%%%

Let $D\subseteq \R^2$ and $\mathbf{p}$ be an interior point of $D$. 
Suppose  $f:D\to\R$ has continuous partial derivatives of first and second order in an open neighborhood of $\mathbf{p}$. Then the  \emph{Hessian form} of $f$ at $\mathbf{p}$ is the binary quadratic form defined by 
$$
Q_{\mathbf{p}}( \mathbf{h})=Q_{\mathbf{p}} (h_1, h_2):=  
f_{xx}\left( \mathbf{p}\right)h_1^2 + 2f_{xy}\left( \mathbf{p}\right)h_1h_2 + f_{yy}\left( \mathbf{p}\right)h_2^2. 
%\frac{\partial^2 f}{\partial x^2}\left( \mathbf{p}\right)h_1^2 
%+  2 \frac{\partial^2 f}{\partial x \partial y}\left( \mathbf{p}\right)  h_1 h_2  +  \frac{\partial^2 f}{\partial %y^2}\left( \mathbf{p}\right)h_2^2 .
$$
%is called the \emph{Hessian form} of $f$ at $\mathbf{p}$. 
%With these notations and hypotheses, we have the following. 
With the hypothesis and notation as above, we have the following. 
\begin{proposition}
\label{Prop:3.2}
If $\nabla f (\mathbf{p}) =\mathbf{0}$ and the Hessian form $Q_{\mathbf{p}}$ of $f$ at $\mathbf{p}$ is indefinite, then $f$ has a strict saddle point at $\mathbf{p}$. 
\end{proposition}

\begin{proof}
%This is proved along similar lines as in most standard 
The basic argument is similar to that used in many texts on Calculus, but we provide a sketch for the sake of completeness.
% (see, e.g., \cite[\S 9.11]{Ap}). 
Assume that $\nabla f (\mathbf{p}) =\mathbf{0}$ and $Q_{\mathbf{p}}$ is indefinite. Then there are nonzero $\mathbf{u}, \mathbf{v}\in \R^2$ such that $Q_{\mathbf{p}}(\mathbf{u})<0$, while $Q_{\mathbf{p}}(\mathbf{v})>0$. 
By the continuity of the second order partials, there is $\delta >0$ 
such that for any $\mathbf{q}\in \R^2$ with $\Vert \mathbf{q} - \mathbf{p} \Vert \le \delta$, we have $\mathbf{q}\in D$ and 
$Q_{\mathbf{q}}(\mathbf{u})<0$, while $Q_{\mathbf{q}}(\mathbf{v})>0$.
%%% More Elaborate Reasoning (footnote?) %%%%%%%%%%%%%%%%%%%%%
%By the continuity of the second order partials, the function 
%$ \mathbf{q} \mapsto Q_{\mathbf{q}}(\mathbf{h})$ is continuous near 
%$ \mathbf{p}$ for any fixed $ \mathbf{h}$. Hence etc.
%%%%%%%%%%%%%%%%%%%%%%%%%%%%%%%%%%%%%%%%%%%%%%%%%%%%%%%%%%%%%%%%% 
%%%%% CHANGE: The following line added
Scaling $\mathbf{u}$ and $\mathbf{v}$ suitably, we may assume that 
$\Vert \mathbf{u} \Vert \le 1$ and $\Vert \mathbf{v} \Vert \le 1$. 
Given any $t\in [-\delta, \delta]$ and $\mathbf{h} \in \R^2$ with 
$\Vert \mathbf{h} \Vert \le 1$, by Taylor's Theorem, there is $\mathbf{q}\in \R^2$ 
on the line joining $\mathbf{p}$ and $\mathbf{p}+ t\mathbf{h}$ such that 
$$
f\left( \mathbf{p} + t\mathbf{h}\right) - f\left(\mathbf{p}\right) 
= \nabla f (\mathbf{p}) \cdot \left(t \mathbf{h}\right) + \frac 12 \, 
Q_{\mathbf{q}}(t\mathbf{h}) = \frac{t^2}{2}\, Q_{\mathbf{q}}(\mathbf{h}).
$$
Thus, if $\gamma_1, \gamma_2: [-\delta,\delta]\to \R^2$ are defined by 
$\gamma_1(t):= \mathbf{p} + t\mathbf{u}$ and 
$\gamma_2(t):= \mathbf{p} + t\mathbf{v}$, then %we see that 
$\gamma_1$ and $\gamma_2$ are regular paths intersecting transversally at $\mathbf{p}$ such  that $f$ has a strict local maximum at $\mathbf{p}$ along $\gamma_1$ and a strict local minimum at $\mathbf{p}$ along $\gamma_2$.
\end{proof}

\begin{remark}
\label{Rem:3.2.1} 
{\rm 
The function $f$ as in Example \ref{Exa:2.3}(ii) has a strict saddle point at $\mathbf{0}$, but its Hessian form at $\mathbf{0}$, being identically zero, is not indefinite. This shows that the converse of  Proposition \ref{Prop:3.2} is not true, in general. We can probe further. Observe that 
our proof of Proposition \ref{Prop:3.2} actually shows that 
%$f$ has a strict saddle point at $\mathbf{p}$. Moreover, it shows that
when the Hessian form is indefinite, the two paths satisfying the 
requirements for a strict saddle point, 
%conditions in Definition \ref{Def:2.2} 
can be chosen as straight line segments. We can, therefore, ask if the `weak converse' is true, that is, if straight line segments suffice to show that a differentiable function has a strict saddle point at $\mathbf{p}$, then whether the Hessian form $Q_{\mathbf{p}}$ is necessarily indefinite? The following example shows that the answer is negative. 
}
\end{remark}

\begin{example} 
\label{Exa:3.2.2}
{\rm
[Dog saddle] The function $f:\R^2 \to \R$ defined by 
$f(x,y):=x^3y-xy^3$ has a saddle point at the origin. To see this, 
it suffices to consider the paths given by $t\mapsto (t, \, -t/2)$ and 
$t\mapsto (t, \, t/2)$. These are straight line segments intersecting transversally at the origin for which the conditions in Definition \ref{Def:2.2} are satisfied. But the Hessian form of $f$ at the origin is identically zero, and hence not indefinite. The graph of $f$ near the origin is shown in Figure \ref{fig:twosaddles} on the right. 
}
\end{example}

In order to apply Proposition \ref{Prop:3.2} to specific examples, it 
is essential to have a useful characterization of the Hessian form 
being indefinite. This is basically a well-known question of Linear Algebra. %(See, for example, \cite{HJ}.) 
(See, e.g., \cite{GL2,HJ}.)  
Again, we include the requisite result and a quick proof for the sake of completeness. 

\begin{lemma}
\label{Lem:3.3}
Let $Q(\mathbf{h}) : = ah_1^2+2bh_1h_2+ch_2^2$ be a binary quadratic form. If $ac-b^2 < 0$,  then $Q$ is indefinite. 
\end{lemma}

\begin{proof}
Observe that (i) if $a\ne 0$, then $Q(1,0)Q(b,-a) = a^2 (ac-b^2)<0$, (ii) if $a=0$ and 
$c\ne 0$, then $Q(0,1)Q(c,-b) = c^2 (ac-b^2) <0$, and (iii) if $a=0$ and $c=0$, then 
$Q(1,1)Q(1,-1) = -4b^2 < 0$. Thus, in any case, $Q$  is indefinite. 
%$$
%$$
%Q(b,-a) = a (ac-b^2) \quad \mbox{and} \quad Q(1,0)=a
%$$
%are nonzero and have opposite signs. Also, if $a=0$ and 
%$c\ne 0$, then 
%$$
%Q(c,-b) = c (ac-b^2) \quad \mbox{and} \quad Q(0,1)=c
%$$
%are nonzero and have opposite signs. Finally, if $a=0$ and $c=0$, then 
%$$
%Q(1,1) = 2b \quad \mbox{and} \quad Q(1,-1)=-2b
%$$
%are nonzero and have opposite signs. Thus, in any case, %we see that 
%$Q$  is indefinite. 
\end{proof}

In the remainder of this section, let $D\subseteq \R^2$ and $\mathbf{p}$ be an interior point of $D$. Further, let 
$f:D\to\R$ be such that $f$ has continuous partial derivatives of first and second order in an open neighborhood 
of $\mathbf{p}$. We define the {\em discriminant} of $f$ at $\mathbf{p}$ to be the real number 
$$
\Delta f (\mathbf{p}) := f_{xx}(\mathbf{p})f_{yy}(\mathbf{p}) - f_{xy}(\mathbf{p})^2 .
$$ 
With the hypothesis and notation as above, we have the following. 

\begin{theorem}[Discriminant Test]
\label{Thm:3.4}
If $\nabla f (\mathbf{p}) =\mathbf{0}$ and $\Delta f (\mathbf{p}) < 0$, then $f$ has a strict saddle point at $\mathbf{p}$.  
\end{theorem}

\begin{proof}
Apply Lemma \ref{Lem:3.3} to the Hessian form %$Q_{\mathbf{p}}$ 
of $f$ at $\mathbf{p}$ and use Proposition \ref{Prop:3.2}.
%Let $H_{\mathbf{p}}$ denote the $n\times n$ real symmetric matrix 
%associated to the Hessian form $Q_{\mathbf{p}}$ of $f$ at $\mathbf{p}$. The hypothesis implies that a principal minor of order 2 of $H_{\mathbf{p}}$ 
%is negative. Hence, by parts (ii) and (iv) of Proposition \ref{Prop:3.3}, 
 %$Q_{\mathbf{p}}$ is indefinite. 
%neither positive semidefinite nor negative semidefinite. Now apply Proposition \ref{Prop:3.2}. 
\end{proof}

%For $1\le i<j\le n$, the real number $f_{x_ix_i}\left( \mathbf{p}\right) f_{x_jx_j}\left( \mathbf{p}\right) - f_{x_ix_j}\left( \mathbf{p}\right)^2$ appearing in Theorem \ref{Thm:3.4} may be called the $(i,j)$-\emph{discriminant} of $f$ at $\mathbf{p}$. Evidently, there are $n(n-1)/2$ such discriminants. According to our Discriminant Test (Theorem \ref{Thm:3.4}), if any one of them is negative, then $f$ has a strict saddle point at $\mathbf{p}$. This is useful in numerous examples. However, as the following example shows, the Discriminant Test can fail but Proposition \ref{Prop:3.2} may still be applicable.
%
%\begin{example}
%\label{exa:3.4.1}
%Consider $f:\R^3\to \R$ defined by 
%%$f(x,y,z):=2xz + 2yz - (x-y)^2 - z^2 - 3x^3 + z^4$. 
%$f(x,y,z):= 4xy + 4xz + 4yz - 3x^2 - 2y^2 - 2z^2$. 
%Clearly $\nabla f\left( \mathbf{0}\right)= \mathbf{0}$ and 
%the Hessian form of $f$ at $\mathbf{0}$ is given by 
%$$
%%Q_{\mathbf{0}}(h,k,\ell)= 4hk+4h\ell + 4k\ell - 2h^2 - 2k^2 - 2\ell^2.
%Q_{\mathbf{0}}(h,k,\ell)= 8hk+8h\ell + 8k\ell - 6h^2 - 4k^2 - 4\ell^2.
%$$
%It is easily seen that the $(1,2)$-, $(1,3)$- and $(2,3)$-discriminants of $f$ are $8$, $0$ and $8$, respectively. 
%%the three possible discriminants of $f$ at $\mathbf{0}$ is zero, and 
%So, the Discriminant Test (Theorem \ref{Thm:3.4}) is not applicable. But $Q_{\mathbf{0}}(1,0,0) = -6 < 0$, while $Q_{\mathbf{0}}(1,1,1) = 10 >0$. Thus $Q_{\mathbf{0}}$ is indefinite and hence by Proposition \ref{Prop:3.2}, $f$ has a strict saddle point at the origin. 
%\end{example}
In fact, as an application of Proposition \ref{Prop:2.4}, we can obtain the following stronger version of Theorem \ref{Thm:3.4}.  % given in the corollary below. 
%test by using the Discriminant Test for local extrema. 

\begin{corollary}
\label{Cor:3.5}
%Let $D\subseteq \R^2$ and $\mathbf{p}$ be an interior point of $D$. 
%Let $f:D\to\R$ be a real-valued function %of two variables $x$ and $y$ 
%with continuous partial derivatives of first and second order in an open neighborhood of $\mathbf{p}$. 
Assume that %the discriminant  
$\Delta f (\mathbf{p})\ne 0$. 
%:= f_{xx}(\mathbf{p})f_{yy}(\mathbf{p}) - f_{xy}(\mathbf{p})^2 $ %of $f$ at $\mathbf{p}$ is nonzero.  
Then 
$$
f \mbox{ has a saddle point at } \mathbf{p} \Longleftrightarrow 
\nabla f (\mathbf{p}) = \mathbf{0} \mbox{ and } \Delta f (\mathbf{p}) < 0.
$$
In particular, when $\Delta f (\mathbf{p})$ is nonzero, $f$ has a saddle point at $\mathbf{p}$ if and only if %$f$ 
it has a strict saddle point at $\mathbf{p}$. 
\end{corollary}

\begin{proof}
If $f$ has a saddle point at $\mathbf{p}$, then it can not have a strict 
local extremum at $\mathbf{p}$. Hence, by the Discriminant test for local extrema of functions of two variables \cite[\S 3.3]{MT}, $\Delta f (\mathbf{p})$ can not be positive. Thus, in view of Proposition \ref{Prop:2.4}, we have 
$\nabla f (\mathbf{p}) = \mathbf{0}$ and $\Delta f (\mathbf{p}) < 0$. The converse follows from Theorem \ref{Thm:3.4}. The last assertion follows from the equivalence just proved and Theorem \ref{Thm:3.4}. 
\end{proof}

\begin{remark}
\label{Rem:3.6}
%The conclusions in 
{\rm 
Example \ref{Exa:2.3}(i) can be treated with the help of 
%alternatively obtained by derived from 
the Discriminant Test. On the other hand, if $f$ is as in Example \ref{Exa:2.3}(ii) or Example \ref{Exa:3.2.2}, then $\Delta f (\mathbf{0}) = 0$, and hence the Discriminant Test is not applicable. 
This shows that the converse of Theorem \ref{Thm:3.4} is not true, in general.
}
\end{remark}

%%%%%%%%%%%%%%%%%%%Older Para %%%%%%%%%%%%%%%%%%%%%%%%%%%%%%%%%%%%%%%%%
%As in Remark \ref{Rem:3.2.1}, the conclusion of Theorem \ref{Thm:3.4} can, in fact, be made stronger to state that $f$ has a strict saddle point at $\mathbf{p}$. Moreover, the desired conditions are achieved 
% via straight line segments intersecting transversally at $\mathbf{p}$. Among the examples we have considered thus far, Theorem \ref{Thm:3.4} is 
%applicable to Example \ref{Exa:2.3}(i). But in the case of Examples 
%\ref{Exa:2.3}(ii) and \ref{Exa:3.2.2}, the discriminant 
%$f_{xx}f_{yy} - f_{xy}^2$ vanishes at the origin and Theorem \ref{Thm:3.4} is not applicable.  However, as we have seen before, the 
%origin is a saddle point in these two examples. This shows that the converse of Theorem \ref{Thm:3.4} is not true, in general. 
%%Further examples are discussed in the next section. 
%%%%%%%%%%%%%%%%%%%%%%%%%%%%%%%%%%%%%%%%%%%%%%%%%%%%%%%%%%%%%%%%%%%%%%

\section{Examples}
\label{sec4}

The aim of this section is to discuss a variety of examples, which not only illustrate our definition of a saddle point but also enable the reader to compare it with the definition usually found in Calculus texts. %We may refer to the latter as the classical definition of saddle point or as saddle points in the classical sense.  
In the latter case, we call it a saddle points in the classical sense.  
%In the process, the pros and cons of our definition of saddle point will %get highlighted. 
%
%To begin with, let us note that although there is a large overlap between the two definitions of saddle points, neither implies the other. For example, 
Note that if $f:\R^2\to \R$ is a constant function and $\mathbf{p}$ is any point of $\R^2$, then  $f$ has a saddle point at $\mathbf{p}$ in our sense but not in the classical sense. On the other hand, if $f:\R^2\to \R$ is defined by $f(x,y):=x^3$, then as we show in Example \ref{Exa:4.1}(iii) below, $f$ has a saddle point at the origin in the classical sense, but not in our sense. However, a strict saddle point in our sense is a saddle point in the classical sense. This is the case when the Discriminant Test (Theorem \ref{Thm:3.4}) is applicable. 
%, one has a saddle point in the classical sense as well as in our sense. 
%More generally, if a function has a strict saddle point, 
%then it has a saddle point in the classical sense. 
%Further contrasting features of the two notions will be apparent from 
%the following examples. 

\begin{examples}
\label{Exa:4.1}
{\rm
(i)  
Consider $f:\R^2 \to \R$ defined by 
$f(x,y):= \min\{|x|, |y|\}$ if $xy\ge 0$ and 
$f(x,y):= - \min\{|x|, |y|\}$ if $xy < 0$.  
%$$
%f(x,y):= \left\{ \begin{array}{rl} \min\{|x|, |y|\} & \mbox{ if } xy\ge 0, \\
%- \min\{|x|, |y|\} & \mbox{ if } xy <  0. \end{array} \right.
%$$
Using the paths $\gamma_1$ and $\gamma_2$ in Example \ref{Exa:2.1}(ii), namely, those given by $t\mapsto (t,-t)$ and $t\mapsto (t,t)$, we see that $f$ has a strict saddle point at $(0,0)$. Thus, a nondifferentiable function can have a saddle point (in our sense). 

(ii) Let $c_1,c_2\in \R$ with $0<c_1<c_2$. Consider $f:\R^2\to \R$ defined by 
$f(x,y):=(y-c_1x^2)(y-c_2x^2)$. Using the paths given by $t\mapsto (t, ct^2)$ and $t\mapsto (0,t)$, where $c\in \R$ satisfies $c_1<c<c_2$, we see that $f$ has a strict saddle point at $(0,0)$. Note that for $\lambda \in \R$, $\lambda\ne 0$, we have $f(t,\lambda t)= t^2(\lambda -c_1t)(\lambda -c_2t) > 0$ for $0< |t|< |\lambda |/c_2$. Also, $f(t,0)= c_1c_2t^4 > 0$ and $f(0,t) = t^2 > 0$ for all 
$t\ne 0$. 
%% CHANGE %It follows that  
Hence
$f$ has a strict local minimum at $(0,0)$ along every straight line through the origin. Thus, in this example, straight line segments alone can not work to show %check 
that $f$ has a saddle point at $(0,0)$, but %although 
a combination of a parabola and a straight line segment does. 

(iii) 
Consider $f:\R^2\to \R$ defined by $f(x,y):=x^m$, where $m$ is an odd positive integer. If $m=1$, then clearly, $f$ has no saddle points (in any sense). Assume now that $m>1$. Then $f$ is 
differentiable and $\nabla f (0,y_0) = (0, 0)$ for any $y_0\in \R$. 
%Also $f$ takes positive and negative values near the origin and hence
Since $f$ takes both positive and negative values in every open neighborhood of $(0,y_0)$, we see that $f$ has a saddle point at $(0,y_0)$ in the classical sense for every $y_0\in \R$. On the other hand, if it had a saddle point (in our sense)  at $(0,y_0)$ for some fixed $y_0\in \R$, then we would find paths $\gamma_i:[a_i,b_i]\to \R^2$ satisfying the conditions of Definition \ref{Def:2.2}. Write $\gamma_i(t):=(x_i(t), y_i(t))$ and let $t_i\in (a_i,b_i)$ be such that  $\gamma_i(t_i) = (0,y_0)$. Then $0$ is a local extremum of $x_i^m$ at $t_i$. But since  $x$ and $x^m$ have the same sign for any $x\in \R$, it follows that each $x_i$ has a local extremum at $t_i$. Consequently, $x_i'(t_i) = 0$ for $i=1, 2$, and therefore $\gamma_1$ and $\gamma_2$ can not intersect transversally at $(0,y_0)$. 
Thus $f$ does not have any saddle point. 
%We remark that this ought not to be surprising if one looks at the graph of $f$ as depicted on the left in  Figure \ref{fig:x2plusy3}. More generally, 
%if $f:\R^2\to \R$ defined by $f(x,y):=x^m$, where $m$ is an odd positive integer, then $\nabla f (0,y_0) = (0, 0)$ but $f$ does not have a saddle point at $(0,y_0)$ for any $y_0\in \R$. This is proved in exactly the same manner as in the case of $x^3$.
}
\end{examples}

Our next set of examples generalize some of the simplest and most natural types of functions of two variables, such as $xy$, $x^2-y^2$, 
$x^3-3xy^2$ and $x^m$, which we have seen earlier. The arguments in the  general case are a bit involved and make good material for starred  exercises in Calculus texts (although we have yet to see them in print), especially for those who
may choose to adopt our definition of saddle point.
% might like to adopt the definition of saddle point that we have proposed.  

\begin{examples}
\label{Exa:4.2}
{\rm 
(i)  Let $m,n\in \N$ (where $\N$ denotes the set of positive integers) and $f:\R^2 \to \R$ be defined by $f(x,y):=x^my^n$. Then
$$
f \mbox{ has a strict saddle point at } (0,0) \Longleftrightarrow 
\mbox{both $m$ and $n$ are odd.}
$$
The implication `$\Longleftarrow$' is easy. Indeed, if  $m$ and $n$ are odd, then $m+n$ is even and $f(t,-t) = -t^{m+n} <0$, while $f(t,t)= t^{m+n} > 0$ for all $t\ne 0$. Thus it suffices to consider the paths $\gamma_1$ and $\gamma_2$ in Example \ref{Exa:2.1}(ii). % do the needful. 
To prove `$\Longrightarrow$', first suppose $m$ and $n$ are both even. Then $f(x,y)\ge 0$ for all $(x,y)\in \R^2$ and hence $f$ can not have a strict local maximum at $(0,0)$ along any path passing through $(0,0)$. So $f$ can not have a strict saddle point at $(0,0)$. Next, suppose $m$ is odd and $n$ is even. For $i=1,2$, let $\gamma_i:[a_i,b_i]\to \R^2$ be regular paths intersecting transversally at $(0,0)$  such that $f$ has a strict local maximum (resp: strict local minimum) along $\gamma_1$ (resp: $\gamma_2$). Write $\gamma_i(t):=(x_i(t), y_i(t))$ and let $t_i\in (a_i,b_i)$ 
be such that $\gamma_i(t_i) = (0,0)$. Then there is  $\delta_1 > 0$ such that 
$$
0<|t-t_1|< \delta_1 \Longrightarrow x_1(t)^my_1(t)^n < 0 \Longrightarrow y_1(t) \ne 0 \mbox{ and } x_1(t) < 0,
$$
where the last implication follows since $n$ is even and $m$ is odd. Thus, $x_1$ has a strict local maximum at $t_1$ and so $x_1'(t_1)=0$. 
Similarly, there is  $\delta_2 > 0$ such that 
$$
0<|t-t_2|< \delta_2 \Longrightarrow x_2(t)^my_2(t)^n > 0 \Longrightarrow y_2(t) \ne 0 \mbox{ and } x_2(t) > 0.
$$
Consequently, $x_2'(t_2)=0=x_1'(t_1)$, which contradicts the assumption that  $\gamma_1$ and $\gamma_2$ intersect transversally at $(0,0)$. The case when $m$ is even and $n$ is odd is similar. Thus, we have shown that if both or one of $m$ and $n$ is even, then $f$ does not have a strict saddle point at $(0,0)$. 

(ii) Let $m,n\in \N$ and $f:\R^2\to \R$ be defined by $f(x,y):=x^m-y^n$. Then
$$
f \mbox{ has a saddle point at } (0,0) \Longleftrightarrow 
\mbox{both $m$ and $n$ are even.}
$$
The implication `$\Longleftarrow$' is again easy. Indeed, if $m$ and $n$ are both even, we have $f(0,t) = -t^{n} < 0$, while $f(t,0)= t^{m} > 0$ for all $t\ne 0$. Thus it suffices to consider the paths $t\mapsto (0,t)$ and $t\mapsto (t,0)$. %  suffice. %do the needful. 
To prove `$\Longrightarrow$', observe that if $m=1$ or $n=1$, then $\nabla f (0,0)\ne (0,0)$, and hence by Proposition \ref{Prop:2.4}, $f$ can not have a saddle point at $(0,0)$.  So we now assume that $m\ge 2$ and $n\ge 2$. 
Suppose $m$ or $n$ is odd.  Note that if $g:\R^2\to \R$ is defined by $g(x,y):=-f(y,x) = x^n-y^m$, then clearly, $f$ has a saddle point at 
$(0,0)$ if and only if $g$ does. Thus, we may assume without loss of generality that $m$ is odd. Let us first consider the case when $n\ge m$. Let $\gamma_i:[a_i,b_i]\to\R^2$ be paths satisfying the conditions in Definition \ref{Def:2.2}.  
Write $\gamma_i(t):=(x_i(t), y_i(t))$ and let $t_i\in (a_i,b_i)$ 
be such that $\gamma_i(t_i) = (0,0)$. Then there is  $\delta_1 > 0$ such that 
$$
|t-t_1|< \delta_1 \Longrightarrow x_1(t)^m \le y_1(t)^n  \Longrightarrow x_1(t) \le y_1(t)^{n/m},
$$
where the last implication follows since $m$ is odd. Since $\gamma_1(t_1) = (0,0)$, we see that $x_1 - y_1^{n/m}$ has a local maximum at $t_1$, and so $x_1'(t_1) - (n/m) y_1(t_1)^{(n-m)/m} y_1'(t_1) = 0$. It follows that $x_1'(t_1) = y_1'(t_1)$ if $n=m$ and $x_1'(t_1) = 0$ if $n>m$. 
Similarly, there is  $\delta_2 > 0$ such that 
$$
|t-t_2|< \delta_2 \Longrightarrow x_2(t)^m \ge y_2(t)^n  \Longrightarrow x_2(t) \ge y_2(t)^{n/m},
$$
and this yields that $x_2'(t_2) = y_2'(t_2)$ if $n=m$ and $x_2'(t_2) = 0$ if $n>m$. Either way, the condition that $\gamma_1$ and $\gamma_2$ intersect transversally is contradicted. 
Next, suppose $m$ is odd and $n < m$. In case $n$ is odd, then considering $g:\R^2\to \R$ defined by $g(x,y):=-f(y,x) = x^n-y^m$, we obtain from the previous case that $g$, and hence $f$, does not have a 
saddle point at $(0,0)$. Thus, let us assume that $n$ is even. Now, as  before, there is $\delta_2>0$ such that 
$$
|t-t_2|< \delta_2 \Longrightarrow x_2(t)^m \ge y_2(t)^n  \ge 0 \Longrightarrow x_2(t) \ge 0,
$$
and this yields $x_2'(t_2)=0$. Consequently, there is $\xi_2:[a_2,b_2]\to \R$ such that $x_2(t) = (t-t_2) \xi_2(t)$ for all $t\in [a_2, b_2]$ and moreover $\xi_2(t) \to 0$ as $t\to 0$. Also, since $y_2$ is differentiable at $t_2$, there is $\eta_2:[a_2,b_2]\to \R$ such that 
$y_2(t) = (t-t_2)\left[ y_2'(t_2)+ \eta_2(t)\right]$ for all $t\in [a_2, b_2]$ and moreover $\eta_2(t) \to 0$ as $t\to 0$. Thus, 
$$
|t-t_2|< \delta_2 \Longrightarrow (t-t_2)^m \xi_2(t)^m \ge 
(t-t_2)^n \left[ y_2'(t_2)+ \eta_2(t)\right]^n ,
$$
and hence
$$
0<|t-t_2|< \delta_2 \Longrightarrow (t-t_2)^{m-n} \xi_2(t)^m \ge 
\left[ y_2'(t_2)+ \eta_2(t)\right]^n .
$$
Since $n$ is even, upon letting $t\to t_2$, we see that 
$0\ge y_2'(t_2)^n = |y_2'(t_2)|^n$, and hence $y_2'(t_2)=0=x_2'(t_2)$. 
So the condition that $\gamma_2$ is regular is contradicted.
%So once again, the condition that $\gamma_1$ and $\gamma_2$ intersect transversally is contradicted.

(iii) Let $m,n\in \N$ and $f:\R^2\to \R$ be defined by $f(x,y):=x^m+y^n$. Then %$\nabla f (0,0)=(0,0)$ but 
$f$ never has a saddle point at $(0,0)$. To see this, note that if $m$ and $n$ are both even, then $f(x,y) > 0$ for all $(x,y)\ne (0,0)$, and so $f$ can not have a local maximum along any path passing through the origin. The remaining cases can be proved by arguments similar to those in (ii) above. 

(iv) [Generalized Monkey Saddle] Let $n\in \N$ and $f:\R^2 \to \R$ be defined by $f(x,y):={\rm Re}\left(x+iy\right)^n$. Note that the surface $z=f(x,y)$ is parametrically given by $x= r\cos \theta$, $y= r\sin \theta$ and $z=r^n\cos n\theta$, where $r\ge 0$ and $-\pi < \theta \le \pi$. 
%For instance, $f(x,y)$ is $x$ or $x^2-y^2$ or $x^3-3x^2y$ when $n$ is %$1$ or $2$ or $3$, respectively. 
If $n=1$, then clearly, $f$ has no saddle points. But $f$ has a strict saddle point at $(0,0)$ if $n\ge 2$.  The case when $n$ is even is easy. In this case one easily sees that it suffices to consider the paths given by $t\mapsto \left(t \cos \left({\pi}/{n}\right), \; t\sin \left({\pi}/{n}\right) \right)$  and  $t\mapsto (t,0)$. 
Next, suppose $n$ is odd and $n>1$. In this case, 
$f$ is negative in the sectors 
%$\pi/2n < \theta < \pi/n$ and $-\pi < \theta < -\pi + \pi/2n$, 
$$
\frac{\pi}{2n} < \theta < \frac{\pi}{n} \quad \mbox{ and } \quad 
-\pi < \theta < -\pi + \frac{\pi}{2n},
$$
whereas $f$ is positive in the sectors 
%$-\pi/2n < \theta < 0$ and $\pi - \pi/n < \theta < \pi - \pi/2n$. 
$$
-\frac{\pi}{2n} < \theta < 0 \quad \mbox{ and } \quad 
\pi -\frac{\pi}{n} < \theta < \pi - \frac{\pi}{2n}.
$$
With this in view, 
%Analyzing the sign of $f$ in the sector $\pi/2n < \theta < \pi/n$, 
we see that if $\gamma_1$ is the parabolic path 
given by 
%$$
%t\mapsto (t,0) \quad \mbox{and} \quad 
%t\mapsto \left(t \cos \frac{\pi}{n}, \; t\sin \frac{\pi}{n} \right),
%$$
%if $n$ is even, and the paths 
$$
t\longmapsto \left(-t \cos \frac{\pi}{2n} + t^2 \sin \frac{\pi}{2n}, \; t\sin \frac{\pi}{2n} + t^2 \cos \frac{\pi}{2n}\right),
$$
then $f(\gamma_1(t))<0$ for $t\ne 0$ with $|t|$ small, while if 
$\gamma_2$ is the parabolic path 
given by 
%\quad \mbox{and} \quad 
$$
t\longmapsto \left(t \cos \frac{\pi}{2n} - t^2 \sin \frac{\pi}{2n}, \; t\sin \frac{\pi}{2n} + t^2 \cos \frac{\pi}{2n}\right), 
$$
then $f(\gamma_2(t)) > 0$ for $t\ne 0$ with $|t|$ small. Also, 
$\gamma_1$ and $\gamma_2$ intersect transversally at $(0,0)$. Thus, $f$ has a 
strict saddle point at $(0,0)$. 

(v) %[Imaginary Monkey Saddle??] 
Let $n\in \N$ and $g:\R^2 \to \R$ be defined by $g(x,y):={\rm Im}\left(x+iy\right)^n$. Note that the surface $z=g(x,y)$  is parametrically given by $x= r\cos \theta$, $y= r\sin \theta$ and $z=r^n\sin n\theta$, where $r\ge 0$ and $-\pi < \theta \le \pi$. By arguments similar to those in (iv) above, it can be proved that $g$ does not have a saddle point at $(0,0)$ if $n=1$, while $g$ has a strict saddle point at $(0,0)$ if $n\ge 2$.
}
\end{examples}

\begin{remark}
\label{Rem:4.3}
{\rm 
Due to the nature of our definition of a saddle point, it is not entirely trivial to show that a specific point is a saddle point of a function. 
The Discriminant Test can help in many cases, but if it fails, then one has to painstakingly construct regular paths with the desired properties. 
To this end, it helps to look at the level curves and to know how the function behaves near the point in question, but still some guessing is needed. Further, when Corollary \ref{Cor:3.5} is not applicable, it becomes even more challenging to show that a point is \emph{not} a saddle point of a given function. For this, we need to logically rule out the existence of regular paths satisfying the properties stated in Definition \ref{Def:2.2}. 
}
\end{remark}

\begin{remark}
\label{Rem:4.4}
{\rm 
As is common in introductory texts on multivariable calculus, we have restricted to %saddle points of 
functions of two variables. 
But it is clear that many of the notions and results discussed here extend readily to $\R^n$ in place of $\R^2$. 
For instance, the notions of paths, regularity, transverse intersections, indefiniteness of a quadratic form, and the Hessian form admit straightforward generalizations. If one uses transversally intersecting regular paths in $\R^n$ to define a saddle point of a function of $n$ variables exactly as in Definition \ref{Def:2.2}, then Proposition \ref{Prop:3.2} continues to hold and Theorem \ref{Thm:3.4} admits an analogue with the condition ``$\Delta f (\mathbf{p}) < 0$'' replaced by 
``$f_{x_ix_i}\left( \mathbf{p}\right)f_{x_jx_j}\left( \mathbf{p}\right) - f_{x_ix_j}\left( \mathbf{p}\right)^2 < 0$ for some $i,j$ with $i\ne j$''. However, the analogue of Proposition \ref{Prop:2.4} for $\R^n$ is not valid if $n>2$.  [Consider, for example, $f:\R^3\to \R$ defined by $f(x_1,x_2,x_3):=x_1x_2+x_3$.] In other words, a saddle point is not automatically a critical point. For this reason, a straightforward analogue of Definition \ref{Def:2.2} is not very satisfactory. A better option may be to define a real-valued function $f$ on an open subset $D$ of $\R^n$ to have a saddle point at $\mathbf{p}\in D$ if there are submanifolds %$\Gamma_1$ and $\Gamma_2$ 
of $D$ whose tangent spaces at $\mathbf{p}$ span $\R^n$ such that $f$ has a local maximum at $\mathbf{p}$ on one and a local minimum at $\mathbf{p}$ on another.  Returning to the case $n=2$, if one wants to consider surfaces more general than those defined by graphs of functions of two variables, another plausible definition for $\mathbf{p}$ to be a  saddle point of a surface $S$ in $\R^3$ could be that there is a plane $P$ passing through $\mathbf{p}$ such that $P\cap S$ is like a ``graph'' for which $\mathbf{p}$ a vertex of degree $\ge 4$. % $z=f(x,y)$
%We believe that the approach proposed here is reasonable and appropriate, especially from a pedagogical point of view,  though we should welcome %However,we will be happy to see expositions of alternative 
%approaches and extensions such as those indicated above. 
}
\end{remark}

%
%
%\section*{Acknowledgments}
%
%We thank an anonymous referee for some helpful remarks. 

\end{document}